# Scoring metrics on separable metric spaces

Kerry M. Soileau
NASA Johnson Space Center
*March 16, 2007*

**Abstract**

We define scoring metrics on separable metric spaces and show that they are always no coarser than the metrics from which they spring.

**1. Introduction.** If one tries to imagine the "simplest" possible metric on a given set $X$, an argument could be made for the trivial metric $\tau_c(x,y) = \begin{cases} c & x \neq y \\ 0 & x = y \end{cases}$, with $c \neq 0$. What is next simplest? We propose the scoring metric. We first envision "delegate functions" $f_n : X \to X$ which associate to each point $x \in X$ a point $f_n(x) \in X$. For any two points $x, y \in X$, for each $n = 1, 2, 3, \cdots$ one may compute the score

$$\tau_{a_n}(f_n(x), f_n(y)) = \begin{cases} a_n & f_n(x) \neq f_n(y) \\ 0 & f_n(x) = f_n(y) \end{cases}, \text{ where } \{a_n\}_{n=1}^{\infty} \text{ is a nonincreasing sequence of}$$

positive real numbers such that $\sum_{i=1}^{\infty} a_i$ converges. These scores are summed to produce the scoring function $\rho(x, y) = \sum_{i=1}^{\infty} \tau_{a_i}(f_i(x), f_i(y))$.

In the following we first show that if $X$ is a separable metric space, it is easy to find functions $f_n : X \to X$ and sequences $\sum_{i=1}^{\infty} a_i$ such that $\sum_{i=1}^{\infty} \tau_{a_i}(f_i(x), f_i(y))$ is a metric on $X$, and that the induced topology is no coarser than the original topology.

**2. Propositions.** Let $X$ be a separable metric space with metric $\sigma(\cdot, \cdot)$. Let $\{r_i\}_{i=1}^{\infty}$ be a countable dense subset of $X$ with $r_i = r_j$ only if $i = j$. Define $f_n : X \to X$ as follows:

$$f_1(x) = r_1,$$

and

$$f_n(x) = \begin{cases} f_{n-1}(x) & \text{if } \sigma(x, f_{n-1}(x)) \leq \sigma(x, r_n) \\ r_n & \text{if } \sigma(x, r_n) < \sigma(x, f_{n-1}(x)) \end{cases}$$

for $n > 1$. Note that

$$\sigma(x, f_n(x)) = \begin{cases} \sigma(x, f_{n-1}(x)) & \text{if } \sigma(x, f_{n-1}(x)) \leq \sigma(x, r_n) \\ \sigma(x, r_n) & \text{if } \sigma(x, r_n) < \sigma(x, f_{n-1}(x)) \end{cases} = \min(\sigma(x, f_{n-1}(x)), \sigma(x, r_n)), \text{ hence}$$

$\sigma(x, f_n(x)) \leq \sigma(x, r_n)$. Fix $\varepsilon > 0$. Since $\{r_i\}_{i=1}^{\infty}$ is dense, we can find $N \geq 1$ such that $\sigma(x, r_N) < \varepsilon$. Since $\sigma(x, f_N(x)) \leq \sigma(x, r_N)$, it follows that $\sigma(x, f_N(x)) < \varepsilon$. By induction on $\sigma(x, f_n(x)) \leq \sigma(x, f_{n-1}(x))$ we get that $\sigma(x, f_n(x)) < \varepsilon$ for all $n \geq N$. Thus $\lim_{n \to \infty} f_n(x) = x$.

Let $\{a_n\}_{n=1}^{\infty}$ be a nonincreasing sequence of positive real numbers such that $\sum_{i=1}^{\infty} a_i$ converges. Then define $\rho: X \times X \to R$ as

$$\rho(x, y) = \sum_{i=1}^{\infty} \tau_{a_i}(f_i(x), f_i(y)).$$

$\rho(\cdot, \cdot)$ is well-defined and finite because it is dominated by $\sum_{i=1}^{\infty} a_i < \infty$.

We claim that $\rho(\cdot, \cdot)$ is a metric over $X$, because of the following:

1. $\rho(x, x) = \sum_{i=1}^{\infty} \tau_{a_i}(f_i(x), f_i(x)) = \sum_{i=1}^{\infty} 0 = 0$.
2. If $x, y \in X$ and $x \neq y$, then given any integer $N$, there exists $n > N$ such that $f_n(x) \neq f_n(y)$, since $\lim_{n \to \infty} f_n(x) = x$ and $\lim_{n \to \infty} f_n(y) = y$. Hence $\rho(x, y) > 0$ whenever $x \neq y$.
3. Fix any three points $x, y, z \in X$ and any positive integer $i$. If $f_i(x) = f_i(z)$ it follows that $\tau_{a_i}(f_i(x), f_i(z)) = 0 \leq \tau_{a_i}(f_i(x), f_i(y)) + \tau_{a_i}(f_i(y), f_i(z))$. If $f_i(x) \neq f_i(z)$, we may infer that either $f_i(x) \neq f_i(y)$ or $f_i(y) \neq f_i(z)$, hence $\tau_{a_i}(f_i(x), f_i(y)) = a_i$ or $\tau_{a_i}(f_i(y), f_i(z)) = a_i$, which in turn implies $\tau_{a_i}(f_i(x), f_i(y)) + \tau_{a_i}(f_i(y), f_i(z)) \geq a_i \geq \tau_{a_i}(f_i(x), f_i(z))$. Hence in either case $\tau_{a_i}(f_i(x), f_i(y)) + \tau_{a_i}(f_i(y), f_i(z)) \geq \tau_{a_i}(f_i(x), f_i(z))$. Next recall that $\rho(x, y) = \sum_{i=1}^{\infty} \tau_{a_i}(f_i(x), f_i(y))$, so

$$\rho(x, y) + \rho(y, z) = \sum_{i=1}^{\infty} \tau_{a_i}(f_i(x), f_i(y)) + \sum_{i=1}^{\infty} \tau_{a_i}(f_i(y), f_i(z))$$

$$= \sum_{i=1}^{\infty} \{\tau_{a_i}(f_i(x), f_i(y)) + \tau_{a_i}(f_i(y), f_i(z))\} \geq \sum_{i=1}^{\infty} \tau_{a_i}(f_i(x), f_i(z)) = \rho(x, z)$$

Thus $\rho(\cdot, \cdot)$ is a metric over $X$.

<u>Proposition 1</u>: Suppose $x_i \xrightarrow{\rho} x$. For any integer $N > 0$, there exists an integer $M > 0$ such that $f_n(x_m) = f_n(x)$ whenever $m \geq M$ and $1 \leq n \leq N$.

Proof: Choose an integer $N > 0$. Since $x_i \xrightarrow{\rho} x$, there exists an integer $M > 0$ such that $\rho(x_m, x) < a_N$ whenever $m \geq M$. Now suppose $f_n(x_m) \neq f_n(x)$ for some $m \geq M$ and some $n$ such that $1 \leq n \leq N$. Then $\rho(x_m, x) = \sum_{\substack{i=1 \\ f_i(x_m) \neq f_i(x)}}^{\infty} a_i \geq a_n \geq a_N > \rho(x_m, x)$ for a contradiction. It then follows that $f_n(x_m) = f_n(x)$ whenever $m \geq M$ and $1 \leq n \leq N$.

Proposition 2: If $a \geq b$ then $\sigma(x, f_a(x)) \leq \sigma(x, r_b)$.

Proof: From the definition, we earlier inferred that
$\sigma(x, f_n(x)) = \min(\sigma(x, f_{n-1}(x)), \sigma(x, r_n))$. This implies $\sigma(x, f_n(x)) \leq \sigma(x, f_{n-1}(x))$, hence by induction $\sigma(x, f_a(x)) \leq \sigma(x, f_b(x))$. Next,
$\sigma(x, f_b(x)) = \min(\sigma(x, f_{b-1}(x)), \sigma(x, r_b)) \leq \sigma(x, r_b)$. Combining these results, we infer $\sigma(x, f_a(x)) \leq \sigma(x, f_b(x)) \leq \sigma(x, r_b)$, and the Proposition is proved.

Proposition 3: If $x, y \in X$ and $x \neq y$, then given any integer $N > 0$, there exists $n > N$ such that $f_n(x) \neq f_n(y)$.

Proof: Note that $f_n(x) \to_\sigma x$ and $f_n(y) \to_\sigma y$ as $n \to \infty$. Now suppose that for some $N > 0$, $f_n(x) = f_n(y)$ for all $n > N$. Then clearly $x = \lim_\sigma f_n(x) = \lim_\sigma f_n(y) = y$, thus $x = y$ for the contradiction.

Proposition 4: For any $\varepsilon > 0$ and $x \in X$, there exists an integer $N > 0$ such that if $f_j(x) = f_j(y)$ for some $j \geq N$ and some $y \in X$, then $\sigma(x, y) < \varepsilon$.

Proof: We prove the contraposition. Choose $\varepsilon > 0$ and $x, y \in X$ such that $\sigma(x, y) \geq \varepsilon$. Because $(X, \sigma)$ is separable, we can find positive integers $m, n$ such that $\sigma(x, r_m) < \frac{\varepsilon}{2}$ and $\sigma(y, r_n) < \frac{\varepsilon}{2}$. Let $N = \max(m, n)$. Suppose $f_j(x) = f_j(y)$ for some $j \geq N$. Then $\sigma(x, f_j(x)) \leq \sigma(x, r_m)$ and $\sigma(y, f_j(y)) \leq \sigma(y, r_n)$, and hence
$\sigma(x, y) \leq \sigma(x, f_j(x)) + \sigma(f_j(x), f_j(y)) + \sigma(y, f_j(y)) < \frac{\varepsilon}{2} + 0 + \frac{\varepsilon}{2} = \varepsilon$ for a contradiction.
Hence $f_j(x) \neq f_j(y)$ for every $j \geq N$.

Theorem: If $x_i \to_\rho x$, then $x_i \to_\sigma x$.

Proof: Pick $\varepsilon > 0$ and $x \in X$. Proposition 4 implies that there exists an integer $N > 0$ such that if $f_j(x) = f_j(y)$ for some $j \geq N$ and some $y \in X$, then $\sigma(x, y) < \varepsilon$. Since $x_i \to_\rho x$, Proposition 1 implies that there exists an integer $M > 0$ such that $f_n(x_m) = f_n(x)$ whenever $m \geq M$ and $1 \leq n \leq N$. In particular, $f_N(x_m) = f_N(x)$ whenever $m \geq M$, hence $\sigma(x, x_m) < \varepsilon$ whenever $m \geq M$, thus $x_i \to_\sigma x$.

Corollary: $(X, \rho)$ is no coarser than $(X, \sigma)$.

Remark: Notice that the Theorem holds regardless of the choice of countable dense subset and of the choice of nonincreasing sequence of positive real numbers with convergent partial sums.